\documentclass[11pt]{amsart}

\newtheorem{thm}{Theorem}

\newtheorem{cor}[thm]{Corollary}

\theoremstyle{remark}
\newtheorem{rem}[thm]{Remark}

\theoremstyle{definition}

\newcommand{\R}{\mathbb{ R}}

\newcommand{\Z}{\mathbb{ Z}}

\title[Bounded cohomology of mapping class groups]{Bounded 
cohomology and non-uniform perfection of mapping class groups}
\author{H.~Endo}
\address{Mathematisches Institut, Universit\"at M\"unchen,
Theresienstr.~39, 80333 M\"unchen, Germany; {\rm current
address:} Department
of Mathematics, Tokyo Institute of Technology, Oh-Okayama, Meguro
152-8551, Tokyo, Japan}
\email{endo@math.titech.ac.jp}
\author{D.~Kotschick}
\address{Mathematisches Institut, Universit\"at M\"unchen,
Theresienstr.~39, 80333 M\"unchen, Germany}
\email{dieter@member.ams.org}
\thanks{Support from the {\sl Deutsche Forschungsgemeinschaft}
is gratefully acknowledged.}
\date{\today; MSC 2000: 57R17, 57R57, 20F12}

\begin{document}

\begin{abstract}
Using the existence of certain symplectic submanifolds
in symplectic $4$-manifolds, we prove an estimate from above
for the number of singular fibers with separating vanishing 
cycles in minimal Lefschetz fibrations over surfaces of positive
genus. This estimate is then used to deduce that mapping class
groups are not uniformly perfect, and that the map from their 
second bounded cohomology to ordinary cohomology is not injective. 
\end{abstract}

\maketitle

\section{Introduction}

It is well-known that the mapping class group of a closed
oriented surface $F$ of genus $\geq 3$ is perfect. 
In this paper we shall prove that it is not uniformly perfect,
meaning that there is no number $N$ such that every element
of the group is a product of at most $N$ commutators. By
a result of Bavard~\cite{B}, this statement implies the 
non-injectivity of the map from the second bounded
cohomology of the mapping class group to its ordinary cohomology. 
This non-injectivity and the non-perfection of the mapping
class groups confirm two conjectures of Morita~\cite{Morita}. 

We shall also prove the above statement about the bounded 
cohomology for some non-perfect subgroups of mapping class groups,
namely for the Torelli groups and for the hyperelliptic mapping 
class groups.

As in~\cite{K}, we use results of Taubes~\cite{Taubes} on the 
Seiberg--Witten invariants of $4$--dimensional symplectic 
manifolds to study the mapping class groups of surfaces. 
We generalise the main argument of~\cite{K} from smooth
surface bundles to Lefschetz fibrations. The result is an
estimate from above for the number of separating singular 
fibers in Lefschetz fibrations over aspherical surfaces, in
terms of the base and fiber genus (and of the number of
nonseparating singular fibers). The existence of such an
estimate implies that mapping class groups are not uniformly 
perfect, even if one does not know the exact shape of the 
estimate. Every concrete estimate, however, quantifies the 
failure of uniformity of perfection;
see Theorem~\ref{t:perfect} below.

\section{The main theorems}

In this section we prove the main technical result about
Lefschetz fibrations and conclude that mapping class groups
are not uniformly perfect. For definitions and further 
information on smooth Lefschetz fibrations we refer to~\cite{GS}.
We assume throughout that Lefschetz fibrations have at most
one critical point in each fiber, and we call the corresponding 
singular fiber separating or nonseparating according to whether
the vanishing cycle is separating or nonseparating.

\begin{thm}\label{t:main}
Let $X$ be a connected smooth closed oriented $4$-manifold and 
$f\colon X\rightarrow B$ a relatively minimal Lefschetz fibration 
with fiber genus $h\geq 2$ and base genus $g\geq 1$ having $s$
separating and $n$ nonseparating singular fibers. Then
\begin{equation}\label{eq:main}
s\leq 6(3h-1)(g-1)+5n \ .
\end{equation}
\end{thm}
\begin{proof}
By a result of Gompf, see~\cite{GS}, $X$ is a symplectic manifold,
and the fibers are symplectic submanifolds. As has been observed
before~\cite{Li,S}, the assumption of relative minimality 
here implies that $X$ is minimal and not ruled, because any 
pseudo-holomorphic sphere in $X$ would have to be contained in a
fiber because $g\geq 1$. Thus Liu's extension~\cite{Liu} of Taubes's 
results~\cite{Taubes} implies $K^2\geq 0$, which we can write as 
\begin{equation}\label{eq:bminus}
b_2^+(X)\geq\frac{1}{5}(b_2^-(X)+4b_1(X)-4) \ .
\end{equation}

Every reducible singular fiber contains a curve of negative 
selfintersection, and all of these are linearly independent
in homology and independent of the class of a smooth fiber,
which has selfintersection zero. Therefore $b_2^-(X)\geq s+1$.
Substituting this in~\eqref{eq:bminus} and using 
$b_1(X)\geq 2g \geq 2$, we obtain
$$
b_2^+(X)\geq 1+\frac{1}{5}s \ .
$$
As the claim~\eqref{eq:main} is trivial for $s=0$, we may
assume $s\geq 1$, and therefore $b_2^+(X)\geq 2$.

The Euler characteristic of $X$ is 
\begin{equation}\label{Euler}
\chi (X) = 4(g-1)(h-1)+s+n \ . 
\end{equation}

We estimate the signature of $X$ using Novikov additivity by 
decomposing the fibration into two pieces. 
Let $D\subset B$ be an embedded $2$-disk containing all the 
critical values of $f$. If $X_1=f^{-1}(D)$, then a result of 
Ozbagci~\cite{O} gives
\begin{equation}\label{O}
\sigma (X_1)\leq n-s \ .
\end{equation} 

Let $X_2$ denote the restriction of $X$ to $B\setminus D$.
As $B\setminus D$ can be decomposed into $2g-1$ pairs of 
pants, and the signature of $X_2$ over each of them is 
given by the Meyer cocycle~\cite{Meyer} and therefore 
bounded by $2h$, we conclude 
\begin{equation}\label{ss}
\sigma (X_2)\leq 2h(2g-1) \ .
\end{equation} 

Combining~\eqref{O} and~\eqref{ss}, we obtain
\begin{equation}\label{sign}
\sigma (X)\leq 2h(2g-1)+n-s \ .
\end{equation}
Passing to finite covers of $B$ and applying~\eqref{sign} to the 
pulled-back fibrations, we finally have
\begin{equation}\label{ssign}
\sigma (X)\leq 2h(2g-2)+n-s \ .
\end{equation}

As $b_2^+(X)\geq 2$, a result of Taubes~\cite{Taubes}
ensures that the canonical class $K$ of $X$ is represented by a 
symplectically embedded surface $\Sigma\subset X$.
It may be disconnected, but because $X$ is minimal, $\Sigma$ has no
spherical component. In the argument below we will tacitly assume
that it is connected. In the general case the same argument works
by summing over the components.

For the genus of $\Sigma$ we have the adjunction formula $g(\Sigma )
=1+K^2=1+2\chi (X)+3\sigma (X)$. Using~\eqref{Euler} and~\eqref{ssign}
we obtain:
\begin{equation}\label{genus}
g(\Sigma )-1\leq 2(10h-4)(g-1) + 5n - s  \ .
\end{equation}
The fibration $f$ induces a smooth map $\Sigma\rightarrow B$ of
degree $d$ equal to the algebraic intersection number of $\Sigma$
with a fiber. This is calculated from the adjunction formula 
applied to a smooth fiber $F$:
\begin{equation}\label{eq:d}
d=\Sigma\cdot F=K\cdot F = 2h-2 \ .
\end{equation}

Now Kneser's inequality $g(\Sigma )-1\geq\vert d\vert (g(B)-1)$
gives:
$$
g(\Sigma )-1\geq 2(h-1)(g-1) \ ,
$$
which together with~\eqref{genus} completes the proof 
of~\eqref{eq:main}.
\end{proof}

The following consequence of Theorem~\ref{t:main} makes 
precise the failure of uniformity of perfection of mapping 
class groups. It also applies to the genus $2$ case, where the 
mapping class group is not perfect. 
\begin{thm}\label{t:perfect}
Let $a$ be a nontrivial separating simple closed curve on a 
surface $F$ of genus $h\geq 2$, and let $t_a$ be the corresponding
Dehn twist. Suppose that $t_a^k$ with $k>0$ can be written as a 
product of $N$ commutators. Then
\begin{equation}\label{eq:bound}
N\geq 1+\frac{k}{6(3h-1)} \ .
\end{equation}
\end{thm}
\begin{proof}
We consider a Lefschetz fibration over the $2$-disk $D$ with precisely
$k$ singular fibers, such that with respect to a basepoint on the 
boundary of $D$ the vanishing cycles of all the singular fibers
can be identified with $a$. Then the monodromy of the fibration around
the boundary of $D$ is $t_a^k$. If this can be expressed as a product
of $N$ commutators, then we can find a smooth surface bundle with fiber
$F$ over a surface of genus $N$ with one boundary component and the 
same restriction to the boundary. Let $X$ be the Lefschetz fibration
over the closed surface $B$ of genus $N$ obtained by gluing together
the two fibrations along their common boundary. 

By construction, no fiber contains a sphere, so $X$ is relatively 
minimal. Thus we can apply Theorem~\ref{t:main} to conclude 
$k\leq 6(3h-1)(N-1)$ as $n=0$ in this case.
\end{proof}

\begin{rem}
Theorem~\ref{t:perfect} implies in particular that no $t_a^k$
equals a single commutator. For $k=1$ this was previously proved 
in~\cite{KO}. The proof there depends on a result of~\cite{S}
whose proof is not correct as written, but can by salvaged 
in the case needed for~\cite{KO}, see the Erratum to~\cite{S}.
\end{rem}

\begin{rem}
It is clear that the number of factors needed to express $t_a^k$ 
as a product of commutators grows at most linearly with $k$. Thus 
Theorem~\ref{t:perfect} settles Problem 2.13 (D) in Kirby's 
list~\cite{Kirby} qualitatively.
\end{rem}

\begin{cor}\label{c:perfect}
Let $\Gamma^k_{h,r}$ be the mapping class group of genus $h$ with
respect to $r$ marked points and $k$ boundary components (fixed 
pointwise). The group $\Gamma^k_{h,r}$ is not uniformly perfect 
for $h\geq 2$.
\end{cor}
\begin{proof}
If $k\geq 1$, then we have a surjective homomorphism 
$$
\Gamma^k_{h,r}\longrightarrow\Gamma^{k-1}_{h,r+1}
$$
given by collapsing a boundary component to a point. We also
have surjective forgetful homomorphisms 
$$
\Gamma^k_{h,r}\longrightarrow\Gamma^k_{h,r-1} \ ,
$$
so it is enough to prove the claim in the case $r=k=0$.
But this case is immediate from Theorem~\ref{t:perfect}. 
\end{proof}

\section{Bounded cohomology and commutator lengths}

In this section we relate Theorems~\ref{t:main} and~\ref{t:perfect}
to the second bounded cohomology.

Let $G$ be a group, and $[G,G]$ its commutator subgroup.
For an element $g\in [G,G]$, the minimal number $c_G(g)$ of 
factors in an expression of $g$ as a product commutators 
is called the commutator length of $g$. The limit
$$
\vert\vert g\vert\vert_G=\lim_{n\rightarrow\infty}\frac{c_G(g^n)}{n}
$$
is called the stable commutator length of $g$. This is related 
to the second bounded cohomology of $G$ by the following result:
\begin{thm}\label{t:B}
{\rm (Bavard~\cite{B})}
The map $H^2_b(G)\rightarrow H^2(G)$ is injective if and only 
if the stable commutator length $\vert\vert g\vert\vert_G$
vanishes identically on $[G,G]$.
\end{thm}

Now the proof of Theorem~\ref{t:perfect} implies that the 
stable commutator length of the Dehn twist along a separating 
simple closed curve is bounded below by $\frac{1}{6(3h-1)}$.
More generally, if $g\in\Gamma^k_{h,r}$ is a product of $s$
separating Dehn twists, not necessarily along the same curve,
then we obtain
\begin{equation}\label{eq:stable}
\vert\vert g\vert\vert_{\Gamma^k_{h,r}}\geq\frac{s}{6(3h-1)} \ .
\end{equation}

Thus Theorem~\ref{t:B} implies:
\begin{cor}\label{c:mapb}
The map $H^2_b(\Gamma^k_{h,r})\rightarrow H^2(\Gamma^k_{h,r})$ 
is not injective for all $h\geq 2$ and $k,r\geq 0$. 
\end{cor}

As Theorem~\ref{t:B} is not limited to perfect groups, we
can also deal with some subgroups of mapping class groups.

\subsection{Hyperelliptic mapping class groups}

Let $\Delta_h$ be the hyperelliptic mapping class group 
of genus $h\geq 2$. It is known that $H^2(\Delta_h,\R )=0$,
cf.~\cite{C,Ka}, so that the statement analogous to 
Corollary~\ref{c:mapb} is just:
\begin{cor}\label{c:hyper}
The bounded cohomology $H^2_b(\Delta_h)$ is non-trivial.
\end{cor}
\begin{proof}
From the presentation of $\Delta_h$ due to Birman-Hilden~\cite{BH},
it follows that the Abelianisation of $\Delta_h$ is a finite
cyclic group of order $4(2h+1)$ if $h$ is odd, and of order 
$2(2h+1)$ if $h$ is even.

Let $a$ be a nontrivial separating simple closed curve on the surface 
of genus $h$ which is invariant under the hyperelliptic involution.
Then $t_a^{4(2h+1)}\in [\Delta_h,\Delta_h]$, and using~\eqref{eq:stable}
we obtain: 
$$
\vert\vert t_a^{4(2h+1)}\vert\vert_{\Delta_h}\geq
\vert\vert t_a^{4(2h+1)}\vert\vert_{\Gamma_h}\geq
\frac{4(2h+1)}{6(3h-1)}>0 \ .
$$
Thus the claim follows from Theorem~\ref{t:B}.
\end{proof}

\subsection{Torelli groups}

Let ${\mathcal T}_h$ be the Torelli group of a closed 
oriented surface of genus $h\geq 2$ defined by
\begin{equation}\label{eq:Torelli}
1\longrightarrow {\mathcal T}_h\longrightarrow\Gamma_{h}
\stackrel{\phi}{\longrightarrow} Sp(2h,\Z )\longrightarrow 1 \ ,
\end{equation}
where $\phi$ denotes the represention of the mapping class 
group on homology.

\begin{cor}\label{c:T}
The map $H^2_b({\mathcal T}_h)\rightarrow H^2({\mathcal T}_h)$ 
is not injective. 
\end{cor}
\begin{proof}
For $h=2$ this follows from the result of Mess~\cite{Mess} that 
${\mathcal T}_2$ is a free group (on infinitely many generators).

Johnson~\cite{J} proved that for $h\geq 3$ the commutator
subgroup of the Torelli group is:
$$
[{\mathcal T}_h,{\mathcal T}_h]={\mathcal T}_h^2\cap {\mathcal K}_h \ ,
$$
where ${\mathcal T}_h^2$ is the subgroup generated by all squares 
of elements of ${\mathcal T}_h$, and ${\mathcal K}_h$ is the 
subgroup generated by the Dehn twists along separating simple closed
curves.

Thus, if $t_a$ is a separating Dehn twist, then
$t_a^2\in [{\mathcal T}_h,{\mathcal T}_h]$.
Using~\eqref{eq:stable}, we have
$$
\vert\vert t_a^2\vert\vert_{{\mathcal T}_h}\geq
\vert\vert t_a^2\vert\vert_{\Gamma_h}\geq
\frac{2}{6(3h-1)}>0 \ ,
$$
and so the claim follows from Theorem~\ref{t:B}.
\end{proof}

\begin{rem}
Here we have used the bound~\eqref{eq:stable} on the stable commutator 
length in $\Gamma_h$ obtained from Theorem~\ref{t:main}, which is 
certainly not optimal for ${\mathcal T}_h$. Examining the proof of
Theorem~\ref{t:main}, we see that for Lefschetz fibrations
whose monodromy is in the Torelli group, we have 
$\sigma (X)\leq n-s$ instead of~\eqref{sign}, so 
that~\eqref{eq:main} is replaced by
\begin{equation}\label{eq:T}
s\leq 6(h-1)(g-1)+5n \ .
\end{equation}
With this we obtain
$$
\vert\vert t_a^2\vert\vert_{{\mathcal T}_h}\geq
\frac{2}{6(h-1)} \ .
$$
\end{rem}

\begin{rem}\label{rem:Burger}
It has been pointed out to us by M.~Burger that the kernel
of $H^2_b(\Gamma_h)\rightarrow H^2(\Gamma_h)$, which is
non-trivial by Corollary~\ref{c:mapb}, injects into the
kernel of $H^2_b({\mathcal T}_h)\rightarrow H^2({\mathcal T}_h)$
under the restriction map. This follows from the exact 
sequence in bounded cohomology associated to an extension
of the form~\eqref{eq:Torelli}, see~\cite{Bou}, together 
with the injectivity in degree $2$ of the map from the bounded to 
the ordinary cohomology of $Sp(2h,\Z )$ proved in~\cite{BM}.
\end{rem}

\bibliographystyle{amsplain}

\bigskip

\end{document}